\documentclass{amsart}

\title[A Self-Dual Hopf Algebra on Double Partially Ordered Sets]{
A Self-Dual Hopf Algebra on Double Partially Ordered Sets }
 \author[C. Malvenuto]{Claudia Malvenuto}
\author[C. Reutenauer]{Christophe Reutenauer}

\address{Claudia Malvenuto\\
Dipartimento  di Informatica \\
Universit\`a La Sapienza \\
Via Salaria 113\\
00198 Roma
Italia} \email{claudia@di.uniroma1.it}

\address{Christophe Reutenauer\\
D\'epartement de Math\'ematiques  UQAM\\
Case Postale 8888 Succ. Centre-ville\\
Montr\'eal (Qu\'ebec) H3C 3P8 Canada
} \email{Reutenauer.Christophe@uqam.ca}

\thanks{This research started
during the first author's stay at LaCIM 2008}


\usepackage{amscd}
\usepackage{verbatim}

\newtheorem{theorem}{Theorem}[section]
\newtheorem{corollary}[theorem]{Corollary}

\newtheorem{prop}{Proposition}[section]
\newtheorem{lemma}{Lemma}[section]

\newfont{\Bbbb}{msbm10 at 12pt}

\newcommand{\BZ}{\mathbb Z}

\newcommand{\BN}{\mathbb N}

\usepackage{graphicx}
\usepackage{pstricks}

\date{\today}
\begin{document}


\begin{abstract}
Let $\mathbf D$ be the set of isomorphism types of finite double
partially ordered sets, that is sets endowed with two partial
orders. On $\BZ\mathbf D$ we define a product and a coproduct,
together with an internal product, that is, degree-preserving. With
these operations $\BZ\mathbf D$ is a Hopf algebra, self-dual with
respect to a scalar product which counts the number of pictures (in
the sense of Zelevinsky) between two double posets. The product and
coproduct correspond respectively to disjoint union of posets and to
a  natural decomposition of a poset into order ideals. Restricting
to special double posets (meaning that the second order is total),
we obtain a notion equivalent to Stanley's labelled posets, and
obtain a sub-Hopf-algebra already considered by Blessenohl and
Schocker. The mapping which maps each double poset onto the sum of
the linear extensions of its first order, identified via its second
(total) order with permutations, is a Hopf algebra homomorphism,
which is isometric and preserves the internal product, onto the Hopf
algebra of permutations, previously considered by the two authors.
Finally, the scalar product between any special double poset and
double posets naturally associated to integer partitions is
described by an extension of the Littlewood-Richardson rule.

\end{abstract}

\maketitle

\section{Introduction} \label{introduction}

We define a (combinatorial) Hopf algebra based on "double posets",
with a scalar product based on "pictures" between double posets, in
analogy to pictures of tableaux as defined by Zelevinsky in
\cite{Zele}. Pictures had been introduced previously by  James and
Peel in \cite{JaPe} p.351-352. Zelevinsky's definition extends
straightforwardly to double posets. The results we prove show that
pictures are fundamentally linked to scalar products, a point of
view already present in Zelevinsky's work, who proved that the
scalar product of two skew Schur functions is equal to the number of
pictures between their shapes. See also \cite{FoGr} and \cite{Leeu}
for the study of pictures between skew shapes.

We call {\em double poset} a set which is endowed with two partial
orders $<_1$ and $<_2$. We consider isomorphism classes of double
posets: on the $\BZ$--module with basis the set of (isomorphism
classes of) double posets, we define combinatorially  a product and
a coproduct, which will make it a graded Hopf algebra, self-dual
with respect to the scalar product $\langle x,y\rangle$ defined as
the number of pictures from $x$ to $y$; in other words,
$$\langle xy,z\rangle=\langle x\otimes y,\delta(z)\rangle.$$
Recall that self-dual Hopf algebras play a great role in representation theory, see \cite{Gei, Zele1}.

When the second order $<_2$ of a double poset is total, one obtains
the  notion which we call special double poset; it is equivalent to
that of labelled poset of Stanley \cite{Stan}, or that of shape of
Blessenohl and Schocker \cite{BlSc}. The corresponding submodule is
then a sub-bialgebra; this bialgebra has already been considered by Blessenohl and Laue, see \cite{BlSc} p.41-42. There is a natural homomorphism into the bialgebra
of permutations of \cite{MaRe}. This mapping is implicit in
Stanley's work (see also \cite{Gess}). The fact that one has a
sub-bialgebra and a homomorphism is already due to \cite{BlSc} (see
also \cite{Mal1} and \cite{HP}). We give further properties of this
homomorphism: it is an isometry, and preserves the internal product.

To each integer partition is naturally associated a special double poset; this construction is described in \cite{Gess}. We describe the scalar product of such a double poset and any special double poset by a rule which extends the Littlewood-Richardson rule.

Note that all the bialgebras in this article are $\BZ$-algebras, are
graded and connected (that is, the 0-component is $\BZ$), hence
these bialgebras are Hopf algebras.

\newcommand{\sh}{\mathbin{\sqcup\!\sqcup}}
\newcommand{\sopra}[2]{\genfrac{}{}{0pt}{}{#1}{#2}}

\section{The bialgebra on double posets}

\subsection{The self-dual bialgebra on double posets}

A {\it double poset} (a notion which is implicit in \cite{MaRe1}) is a triple $(E,<_1,<_2)$, where $E$ is a {\it
finite} set, and $<_1$ and $<_2$ are two partial orders on $E$. When
no confusion arises, we denote  $(E,<_1,<_2)$ simply by $E$. We call
$<_1$ the {\it first order} of $E$ and $<_2$ {\it the second order}
of $E$.

As expected, we say that two double posets $(E,<_1,<_2)$ and
$(F,<_1,<_2)$ are {\it isomorphic} if there exists a bijection
$\phi: E\rightarrow F$ which is an isomorphism from the partial
order $(E,<_1)$ to $(F,<_1)$ and from $(E,<_2)$ to $(F,<_2)$, i.e.
$$\forall x,y \in E: \
x <_i y  \mbox{ in } E \Leftrightarrow \phi(x)<_{i}\phi(y) \mbox{ in } F, \mbox{ for }i=1,2.$$

Rather than on double posets, we want to work on isomorphism classes
of double posets: to avoid too much notation, we simply say double
poset, meaning its isomorphism class.

Let ${\mathbf D}$ denote the set of double posets. We define several
combinatorial operations on this set, which will serve to define the
bialgebra structure on $\BZ{\mathbf D}$, the set of $\BZ$-linear
combinations of double posets.

If $E$ and $F$ are two double posets, their {\it composition}, denoted $EF$,
is the double poset \mbox{$(E\cup F,<_1,<_2)$},  where the union is disjoint and where
\begin{itemize}
\item the first order $<_1$  of $EF$ is the extension to $E\cup F$ of the
first orders $<_1$ of $E$ and $F$, and no element of $E$ is
comparable to any element of $F$;
\item the second order $<_2$ of $EF$ is
the extension to $E\cup F$ of the second orders $<_2$ of $E$ and
$F$, together with the relations $e<_2 f$ for any $e\in E, f\in
F$.
\end{itemize}
The product on $\BZ{\mathbf D}$ is obtained by extending linearly the composition on
$\mathbf D$.

Recall  that an {\it inferior ideal} of a poset $(E,<)$ is a subset $I\subseteq E$ such that if
$y\in I$ and $x< y$, then $x\in I$. A {\it superior ideal} of  $E$ is a subset
$S\subseteq E$ such that if $x\in J$ and $x< y$, then $y\in S$.
Clearly, the complement of an inferior ideal is a superior ideal and conversely.
A {\it decomposition} of a poset $(E,<)$ is a couple $(I,S)$
where $I$ is an inferior ideal and $S$ its complement.

We call {\it decomposition}  of a double poset $(E,<_1,<_2)$ a pair
$$((I,<_1,<_2),(S,<_1,<_2)),$$ where $(I,S)$ is a decomposition of the poset $(E,<_1)$, and where the
first and second orders $<_1,<_2$ for $I$ and
$S$ are obtained by restricting the orders $<_1,<_2$ of the double poset
$E$.

Now let $\delta: \BZ{\mathbf D}\longrightarrow \BZ{\mathbf D}\otimes
\BZ{\mathbf D}$ be the linear map defined on ${\mathbf D} $ by
\begin{equation}\label{defcopr}
\delta((E,<_1,<_2))=\sum_{}(I,<_1,<_2)\otimes (S,<_1,<_2),
\end{equation}
where the sum is extended to all decompositions of $(E,<_1,<_2)$.

A {\it picture} between double posets $(E,<_1,<_2)$ and
$(F,<_1,<_2)$ is a bijection $\phi:E\rightarrow F$ such that:
\begin{itemize}
\item $e<_1 e' \Rightarrow \phi(e)<_2 \phi(e')$ and
\item $f<_1 f' \Rightarrow \phi^{-1}(f)<_2 \phi^{-1}(f').$
\end{itemize}
In other words, a picture is a bijection $\phi$ of $E$ to $F$ which
is increasing from the first order of $E$ to the second order of $F$
and such that its inverse $\phi^{-1}$ is increasing from the first
order of $F$ to the second order of $E$. We define a pairing
$\langle ,\rangle: \BZ{\mathbf D}\times \BZ{\mathbf D}\rightarrow
\BZ$ for any double posets $E,F$ by:
$$\langle E,F\rangle=|\{\alpha: E\rightarrow F, \alpha \mbox{ is a picture}\}|,$$
and extend it bilinearly to obtain a scalar product on $\BZ{\mathbf
D}$, which we call {\it Zelevinsky scalar product}.

\begin{theorem} $\BZ{\mathbf D}$ is a graded
self-dual Hopf algebra.
\end{theorem}

Note that a similar bialgebra structure has been defined on posets
(not double posets) by Schmitt in \cite{Schm} p. 27-28; see also
\cite{ABSo} Example 2.3. Formally it means that the mapping of
$\BZ{\mathbf D}$ into the bialgebra of Schmitt, which sends a double
poset onto the poset with only the first order $<_1$, is a bialgebra
homomorphism.

Self-duality means that for any double posets $E,F,G$,
$$\langle EF,G\rangle=\langle E\otimes G, \delta G\rangle.$$

\noindent{\bf Proof.} We omit the easy verification of the
associativity of the product and of the coassociativity of the
coproduct. Similarly for the homogeneity of both, where the degree
of a double poset is the number of its elements. In order to show
that $\BZ{\mathbf D}$ is a bialgebra, we show that the coproduct
$\delta$ is a homomorphism for the product. This is rather a
tautology, since isomorphic double posets are identified. It amounts
to show that there is a bijection between the set of decompositions
of the double poset $EF$ and the set of pairs $(E_iF_i,E_sF_s)$,
where $(E_i,E_s)$ is a decomposition of $E$ and $(F_i,F_s)$ is a
decomposition of $F$. The bijection is the natural one: take a
decomposition $(I,S)$ of $EF$; then $I$ is an inferior ideal of
$(E\cup F,<_1)$; hence, $I\cap E$ is an inferior ideal of $(E,<_1)$,
$I\cap F$ is an inferior ideal of $(E,<_1)$ and $I=(I\cap E)(I\cap
F)$. Similarly, $S\cap E$ is a superior ideal of $(E,<_1)$, $S\cap
F$ is a superior ideal of $(F,<_1)$ and $S=(S\cap E)(S\cap F)$.
Hence the mapping is the identity mapping (modulo double posets
isomorphisms):
$$(I,S)\mapsto \big((I\cap E)(I\cap F), (S\cap E)(S\cap F)\big).$$
To see that it is a bijection, note that if $(E_i,E_s)$ is a
decomposition of $E$ and $(F_i,F_s)$ is a decomposition of $F$, then
$E_i\cup F_i$ is an inferior ideal $I$, and $E_s\cup F_s$ the
complementary superior ideal $S$ of $(E\cup F, <_1)$; moreover,
$I=E_iF_i$ and $S=E_sF_s$.

We prove now self-duality. In view of the identity stated before the
proof, this amounts to give, for any double posets $E,F,G$, a
bijection between pictures from $EF$ to $G$ and 4-tuples
$(\phi,\psi,I,S)$, where $I$ is an inferior ideal of $G$, $S$ the
complementary superior ideal, $\phi$ a picture of $E$ onto $I$ and
$\psi$ a picture of $F$ onto $S$. So, let $\alpha$ be a picture from
$EF$ onto $G$. Define $I=\alpha(E)$, $S=\alpha(F)$ and the
bijections $\phi: E \rightarrow I, \psi:F\rightarrow S$ obtained by
restriction of $\alpha$ to $E$ and $F$. We verify first that $I$ is
an inferior ideal of G; take $g,g'$ in G with $g<_1 g'$ and $g'\in
I$, hence $\alpha^{-1} (g')\in E$. Then, $\alpha$ being a picture,
we have $\alpha^{-1}(g)<_2 \alpha^{-1}(g')$. Now, if we had $g\notin
I$, then $g\in S$, hence $\alpha^{-1}(g) \in F$, hence
$\alpha^{-1}(g')< _2 \alpha^{-1}(g)$, by definition of the second
ordre of $EF$: contradiction. Similarly, $S$ is a superior ideal of
$G$. Now, the restriction of a picture is a picture, so $\phi$ and
$\psi$ are pictures.

Conversely, given a 4-tuple as above, we glue together the two
bijections $\phi$ and $\psi$ and obtain a bijection $\alpha: E
\rightarrow G$. Since in $EF$, elements of $E$ and elements of $F$
are $<_1$-incomparable, the fact that $\alpha$ is increasing from
$(EF,<_1)$ onto $(G,<_2)$ follows from the similar property for
$\phi$ and $\psi$. Now let $g,g'\in G$ with $g<_1 g'$; if they are
both in $I$ or both in $S$, then $\alpha^{-1}(g)<_2\alpha^{-1}(g')$,
by the similar property for $\phi$ and $\psi$; otherwise, we have
$g\in I$ and $g'\in S$, since $I,S$ are ideals. Then
$\alpha^{-1}(g)\in E$ and $\alpha^{-1}(g')\in F$, and consequently,
$\alpha^{-1}(g)<_2\alpha^{-1}(g')$, by the definition of the second
order of $EF$. Thus $\alpha$ is a picture.

The counity map $\epsilon :\BZ{\mathbf D} \longrightarrow \BZ $ maps
the empty double poset on $1$, and all other doubles posets onto
$0$. It is a morphism for the product. By a well--known fact, a
graded connected bialgebra is a Hopf algebra. \hfill $\Box$

\subsection{A homomorphism into quasi-symmetric functions}

Let $\pi=(E,<_1,<_2)$ be a double poset. Similarly to \cite{Stan}
and \cite{Gess}, we call {\it $\pi$-partition}, a function $x$ from
$E$ into a totally ordered set $X$, such that:
\begin{itemize}
\item $e <_1 e'$  implies $x(e)\leq x(e')$;
\item $e <_1 e'$ and $e\geq_2 e'$ implies $x(e)<x(e')$.
\end{itemize}

Note that if $x$ is injective, the first condition suffices.

Now suppose that $X$ is an infinite totally ordered set of commuting
variables. Then the generating quasi-symmetric function of $\pi$ is
the sum, over all $\pi$-partitions, of the monomials $\prod_{e\in E}
x(e)$. We denote it $\Gamma (\pi)$. By extending linearly $\Gamma$
to $\BZ{\mathbf D}$, we obtain a linear mapping into the algebra of
quasi-symmetric function. For quasi-symmetric functions, see
\cite{Stan1} 7.19. They form a bialgebra denoted $\mathbf {QSym}$,
see e.g. \cite{MaRe}.

\begin{theorem} $\Gamma:\BZ{\mathbf D}\rightarrow \mathbf {QSym}$ is a homomorphism of
bialgebras.
\end{theorem}

This result is implicit in
\cite{Mal}, Proposition 4.6 and Theor\`eme 4.16. We therefore omit the proof, recalling only the
definition of the coproduct of the bialgebra of quasi-symmetric
functions (\cite{Gess} p.300). Take a second
infinite totally ordered set of commuting variables Y. Now order
$X\cup Y$ by $x<y$ for any $x\in X$  and $y\in Y$.  Then a
quasi-symmetric function $F$ on $X$ defines a unique quasi-symmetric
function on $X\cup  Y$, which may be rewritten as a finite linear
combination $\sum_i F_i(X)G_i(Y)$, for some quasi-symmetric
functions $F_i$, $G_i$. Then the coproduct $\delta(F)$ is defined as
$\sum_i  F_i \otimes G_i$.

Note also that a bialgebra constructed on posets has been considered
in \cite{Ehre} and \cite{BeSo}, together with a bialgebra
homomorphism into quasi-symmetric functions. There seems to be no
evident link between their construction and ours, although the
coproduct is essentially the same, while the product rests on
product of posets.

\subsection{An internal product}

Let $(E,<_1,<_2)$ and $(F,<_1,<_2)$ be two double posets. Let $\phi:
(E,<_1)\rightarrow (F,<_2)$ be a bijection  and denote its graph by
$$E\times_\phi F=\{(e,f): \phi(e)=f\}.$$
This set becomes a double poset by: $(e,f)<_1(e',f')$ if and only if
$f<_1f'$; $(e,f)<_2(e',f')$ if and only if $e<_2e'$. In other words,
denoting by $p_1, p_2$ the first and second projections, $p_1$ is an
order isomorphism $(E\times_\phi F,<_2)\rightarrow (E,<_2)$ and
$p_2$ is an order isomorphism $(E\times_\phi F,<_1)\rightarrow
(F,<_1)$. Note that the inverse isomorphisms are
$p_1^{-1}=(id,\phi)$ and $p_2^{-1}=(\phi^{-1},id)$.

Define the {\it internal product} of $(E,<_1,<_2)$ and $(F,<_1,<_2)$
as the sum of the double posets $E\times_\phi F$ for all increasing
bijections  $\phi: (E,<_1)\rightarrow (F,<_2)$. It is denoted
$E\circ F$.

Note that this product has been chosen, among several symmetrical
ones, so that the following holds: let $\sigma$ be a permutation in
$S_n$ and denote by $P_\sigma$ the double poset with underlying set
$\{1,...,n\}$, with $<_2$ as the natural order of this set, and
$<_1$ the order defined by
$\sigma(1)<_1\sigma(2)<_1...<_1\sigma(n)$. Then given two
permutations $\tau$ and $\sigma$, one has
$$P_\sigma\circ P_\tau = P_{\sigma\circ\tau},$$
as the reader may easily verify.


The internal product is compatible with the Zelevinsky scalar product, as follows.

\begin{prop}
Let $E,F,G$ be double posets. Then
$\langle E\circ F,G \rangle=\langle E,F\circ G\rangle$.
\end{prop}

The proposition immediately follows from the following lemma.

\begin{lemma}
Let $(E,<_1,<_2)$, $(F,<_1,<_2)$, $(G,<_1,<_2)$ be three double
posets. There is a natural bijection between

(i) the set of pairs  $(\phi,\alpha)$, where $\phi$ is an increasing
bijection $(E,<_1)\rightarrow (F,<_2)$ and and $\alpha$ is a picture
from $E\times_\phi F$ into $G$;

(ii) the set of pairs $(\psi,\beta)$, where $\psi$ is an increasing
bijection $(F,<_1)\rightarrow (G,<_2)$ and  $\beta$ a picture from
$E$ into $F\times_\psi G$.
\end{lemma}

\noindent{\bf Proof.} We show that the bijection is defined by
$\psi=\alpha\circ(\phi^{-1},id)$ and $\beta=(id,\psi)\circ\phi$, the
inverse bijection being defined by $\phi=p_1\circ \beta$, with $p_1$
the projection $F\times G\rightarrow  F$ and $\alpha=\psi\circ p_2$
with $p_2$ the projection $E\times F\rightarrow F$.

1. Let $(\phi,\alpha)$ be as in (i) and define
$\psi=\alpha\circ(\phi^{-1},id)$ and $\beta=(id,\psi)\circ\phi$. Now
notice that  $(\phi^{-1},id)$ is a mapping $F\rightarrow
E\times_\phi F$ and by definition of the first order of
$E\times_\phi F$, it is increasing for the first orders on $F$ and
$E\times_\phi F$. Since $\alpha$ is increasing $(E\times_\phi
F,<_1)\rightarrow (G,<_2)$, we see that $\psi$ is increasing
$(F,<_1)\rightarrow (G,<_2)$. Now $\beta$ maps bijectively $E$ into
$F\times_\phi G$, as desired, and we verify that it is a picture.
Note that $(id,\psi):F\rightarrow F\times_\psi G$ is increasing for
the second orders, by the definition of the latter on $F\times_\psi
G$. Hence $\beta$ is increasing $(E,<_1)\rightarrow (F\times_\psi
G,<_2)$. Moreover, let $(f,g), (f',g')\in F\times_\psi G$ with
$(f,g)<_1(f',g')$; then $\psi(f)=g$, $\psi(f')=g'$ and  $g<_1g'$.
Thus $\psi(f)<_1\psi(f')$. We have for some $e,e'\in E$,
$\beta(e)=(f,g)$ and $\beta(e')=(f',g')$. Note that the definition
of $\psi$ implies that $\alpha^{-1}\circ\psi=(\phi^{-1},id)$, thus
$\alpha^{-1}(\psi(f)=(\phi^{-1}(f),f)$. Since $\alpha$ is a picture,
$\alpha^{-1}$ is increasing $(G,<_1)\rightarrow (E\times_\phi
F,<_2)$, so that $\alpha^{-1}(\psi(f))<_2\alpha^{-1}(\psi(f'))$ and
therefore by definition of the second order on $E\times_\phi F$, we
have $\phi^{-1}(f)<_2\phi^{-1}(f')$; since
$\beta(e)=(\phi(e),\psi(\phi(e))$, we obtain $f=\phi(e)$ and
$e=\phi^{-1}(f)$. Finally $e<_2e'$ showing that $\beta$ is a
picture.

Define now $\phi'=p_1\circ\beta$ and $\alpha'=\psi\circ p_2$. We must show that
$\phi'=\phi$ and $\alpha'=\alpha$. We have $\phi'=p_1\circ (id,\psi)\circ\phi$,
which is equal to $\phi$ since $p_1\circ(id,\psi)$ is the identity of $F$.
Moreover, $\alpha'=\alpha\circ (\phi^{-1},id)\circ p_2$, and we are done,
since $(\phi^{-1},id)\circ p_2$ is the identity of $E\times_\phi F$.

2. Let $(\psi,\beta)$ be as in (ii), and define $\phi=p_1\circ \beta$
and $\alpha=\psi\circ p_2$. Now,  $\beta$ is an increasing bijection
$(E,<_1)\rightarrow (F\times_\psi G,<_2)$ and $p_1$ is an increasing bijection
$(F\times_\psi G,<_2)\rightarrow (F,<_2)$ by definition of the second order on
$F\times_\psi G$. Thus $\phi$ is an increasing bijection $(E,<_1)\rightarrow(F,<_2)$.
Moreover, $p_2$ is an increasing bijection $(E\times_\phi F,<_1)\rightarrow (F,<_1)$
by definition of the first order of $E\times_\phi F$, and $\psi$
is an increasing bijection $(F,<_1)\rightarrow (G,<_2)$.
Thus $\alpha$ is an increasing bijection:
$$\alpha:(E\times_\phi F,<_1)\rightarrow(G,<_2).$$ We show that $\alpha^{-1}$
is also increasing $(G,<_1)\rightarrow(E\times_\phi F,<_2)$. Indeed, let $g,g'\in G$
with $g<_1g'$. Then $g=\alpha(e,f)$, $g'=\alpha(e',f')$ with $\phi(e)=f$,
$\phi(e')=f'$ and we must show that $(e,f)<_2(e',f')$, that is, $e<_2e'$.
Since $\beta$ is a picture, $\beta^{-1}$ is increasing $(F\times_\psi G,<_1)\rightarrow(E,<_2)$.
We have $g=\alpha(e,f)=\psi(f)$ and similarly $g'=\psi(f')$. Hence $(f,g), (f',g')\in F\times_\psi G$
and $(f,g)<_1(f',g')$ since $g<_1g'$. Thus $\beta^{-1}(f,g)<_2\beta^{-1}(f',g')$.
Now $f=\phi(e)=p_1(\beta(e))\Rightarrow\beta(e)=(f,\psi(f))=(f,g)\Rightarrow e=\beta^{-1}(f,g)$.
Similarly, $e'=\beta^{-1}(f',g')$. It follows that  $e<_2e'$.

Define now $\psi'=\alpha\circ(\phi^{-1},id)$ and $\beta'=(id, \psi)\circ \phi$.
We need to show that $\psi'=\psi$ and $\beta'=\beta'$.
We have $\psi'=\psi\circ p_2\circ(\phi^{-1},id)$ which is clearly equal to $\psi$.
Moreover, $\beta'=(id,\psi)\circ p_1\circ\beta$ and we are done since
$(id,\psi)\circ p_1$ is the identity of $F\times_\psi G$.
\hfill $\Box$

\section{The sub-bialgebra of special double posets}

\subsection{Special double posets}

We call a double poset {\it special} if its second order is total.
Since we identify isomorphic double posets, a special double poset
is nothing else than a labelled poset in the sense of \cite{Stan}.
Given a labelled poset, the labelling (which is a bijection of the
poset into $\{1,\ldots,n\}$) defines a second order, which is total,
on the poset. We denote by $\BZ\mathbf {DS}$ the submodule of
$\BZ{\mathbf D}$ spanned by the special double posets.

A {\it linear extension} of a special double poset $\pi=(E,<_1,<_2)$
is a total order on $E$ which extends the first order $<_1$ of $E$.
We may identify a total order on $E$ with the word obtained by
listing the elements increasingly for this order: let $e_1e_2\ldots
e_n$ be this word, with $\mid E\mid =n$. Let moreover $\omega$ be
the {\it labelling} of $\pi$, that is the unique order isomorphism
from $(E,<_2)$ onto $\{1,\ldots,n\}$. Then we identify the linear
extension with the permutation $\sigma=\omega(e_1)...\omega(e_n)$ in
$S_n$; in other words $\sigma(i)=\omega(e_i)$ and the mapping
$\sigma^{-1}\circ\omega$ is an increasing bijection
$(E,<_1)\rightarrow \{1,\ldots,n\}$, since $\sigma^{-1}\circ\omega(e_i)=i$. In this way, a linear extension
of $\pi$ is a permutation $\sigma\in S_n$ such that
$\sigma^{-1}\circ\omega$ is an increasing bijection
$(E,<_1)\rightarrow \{1,\ldots,n\}$.

In \cite{MaRe} a bialgebra structure on $\BZ  S =\oplus_{n\in \BN}
\BZ S_n$ has been constructed. We recall it briefly. Recall that,
for any word $w$ of length $n$ on a totally ordered alphabet, the
{\it standard permutation} of $w$, denoted by $st(w)$, is the
permutation which is obtained by giving the numbers $1,\ldots,n$ to
the positions of the letters in $w$, starting with the smallest
letter from left to right, then the second smallest, and so on. For
example, if  $w=4\ 3\ 2\ 4\ 1\ 3\ 4\ 4\ 2\ 3\ 3$, then $st(w)=8\ 4\
2\ 9\ 1\ 5\ 10\ 11\ 3\ 6\ 7$. The product $\ast$ for two
permutations $\sigma\in S_n$ and $\tau\in S_p$ is defined as the sum of
permutations in $S_{n+p}$ which are in the {\it shifted shuffle
product}  of the words $\sigma$ and $\tau$, that is the shuffle
product of $\sigma$ and $\bar \tau$, where the latter word is obtained
from $\tau$ by replacing in it each digit  $j$ by $j+n$. The
coproduct $\delta$ on $\BZ  S$ is defined on a permutation
$\sigma\in S_n$ by: $\delta(\sigma)$ is the sum, over all
factorizations (as concatenation) $\sigma=uv$ of the word $\sigma$,
of $st(u)\otimes st(v)$, where $st$ denotes standardization of a
word. See \cite{MaRe} for these definitions.

On $\BZ  S$ put the {\it J\"{o}llenbeck scalar product}, defined by
\[
(\sigma,\tau)=
\left\{
\begin{array}{cl}
 1 &  \mbox{ if }   \sigma=\tau^{-1}\\
 0 &     \mbox{ otherwise. }
\end{array}
\right.
\]
It turns $\BZ  S$ into an self-dual bialgebra, see \cite{BlSc} 5.14
(where an isomorphic algebra is considered, obtained by replacing each permutation by its inverse).

Define the linear mapping $L$ from $\BZ\mathbf {DS}$  into $\BZ  S$
by sending each special double poset on the sum of its linear
extensions. The following result is already in \cite{BlSc}, see
4.18, 5.5 and 5.10. Note that what we call here special double poset
is called {\em shape} by Blessenohl and Schocker. Also, what we call
composition is called {\em semi-direct product} by them. For sake of
completeness, we give here an alternative  proof of the coalgebra
property; we construct a bijection between the appropriate sets of
decompositions. The proof in \cite{BlSc} (different from ours) uses an argument, due to
Gessel \cite{Gess}, which is inductive on the  number of edges in
the Hasse diagramm of the poset.

\begin{theorem} $\BZ\mathbf {DS}$ is a sub-bialgebra of
$\BZ{\mathbf D}$ and $L:\BZ\mathbf {DS}\rightarrow\BZ  S$ is homomorphism of bialgebras.
\end{theorem}

\noindent{\bf Proof.} It is straightforward to see that the
composition of two special double posets is special, and that a
decomposition of a special double poset is a pair of special double
posets (the class of special double posets is a {\it hereditary
family} in the sense of Schmitt \cite{Schm}: it is closed under
taking disjoint unions and ideals). Thus, $\BZ\mathbf {DS}$ is a
sub-bialgebra of $\BZ{\mathbf D}$. Moreover, the set of linear
extensions of the composition of $\pi$ and $\pi'$ is classically the
shifted shuffle product of the set of linear extensions of $\pi$ by
that of $\pi'$. Hence $L$ is a homomorphism of algebras.

The fact that it is also a homomorphism of coalgebras is proved as
follows. Let $\pi=(E,<_1,<_2)$ be a special double poset. Then
$$\delta\circ L(\pi)=\sum_{\sigma,u,v}st(u)\otimes st(v),$$ where
the summation is over all  triples $(\sigma,u,v)$, with $\sigma$ a
linear extension of $(E,<_1)$ and where $\sigma$ is the
concatenation $uv$; moreover, $$(L\otimes L)\circ\delta
(\pi)=\sum_{I,S,\alpha,\beta}\alpha\otimes\beta,$$ where the
summation is over all quadruples $(I,S,\alpha,\beta)$ with $I$ an
inferior ideal of $(E,<_1)$ and $S$  its complementary superior
ideal, and $\alpha, \beta$ are respectively linear extensions of $I,
S$ for the induced order $<_1$. We show that there is a bijection
between the set of such triples and quadruples. To simplify,  take
$\pi=(E,<_1,<_2)$ with $E=\{1,...,n\}$ and $<_2=<$ the natural total
order on $E$. Then the labelling $\omega$ of $\pi$ is the identity
mapping. We show that $$(\sigma,u,v)\mapsto  (I,S,st(u),st(v)),$$
with $I$ the set of naturals appearing in $u$ and $S$ the set of
naturals appearing in  $v$, is the desired bijection. Note first
that, since $\sigma$ is a linear extension of $(E,<_1)$ and
$\sigma=uv$, then $I$ and $S$ as defined are a lower and a superior
ideal of $(E,<_1)$; moreover,  $st(u)$ and $st(v)$ are linear
extensions of $(I,<_1)$ and $(S,<_1)$. This mapping is injective,
since any permutation $\sigma=uv$ is determined by $st(u)$, $st(v)$
and the sets of digits in $u$ and $v$. We show that it is also
surjective: let $(I,S,\alpha,\beta)$ a quadruple as above, and
define uniquely $\sigma=uv$ with $st(u)=\alpha$, $st(v)=\sigma$ and
$I,S$ the set of digits in $u,v$. All we have to show is that
$\sigma$ is a linear extension of $(E,<_1)$. That is: if
$e=\sigma(j)$ and $e'=\sigma(k)$, with $e<_1e'$, then $j<k$. Since
$\sigma$ is the concatenation of $u$  and $v$, and since $st(u)$,
$st(v)$ are linear extensions of $I,S$  for the order $<_1$, this is
clear if $j,k$ are both digits in $\{1,...,i\}$ or
$\{i+1,...,i+s\}$, with $i,s$ the cardinality of  $I,S$; also, if
$j$ is in the first set and $k$  in the second. Suppose by
contradiction that $j$ is in the second set and $k$ in the  first;
then $e\in S$ and $e'\in I$, contradicting the ideal property.
 \hfill $\Box$

\bigskip

The homomorphism $L$ has two other properties.

\begin{theorem} The homomorphism $L:\BZ\mathbf {DS}\rightarrow\BZ  S$
preserves the Zelevinsky scalar product and the internal product.
\end{theorem}

Recall that the {\it internal product} of $\BZ  S$ is simply the
product which extends the product on permutations. The next lemma
extends known results on classical pictures between skew shapes, cf. \cite{GR},
\cite{BlSc} Remark 13.6.

\begin{lemma}
Let $\pi,\pi'$ be special double posets. There is a natural
bijection between pictures from $\pi$ into  $\pi'$ and linear
extensions of $\pi$ whose inverse is a linear extension of $\pi'$.
\end{lemma}
\noindent{\bf Proof.} Let $\phi$ be a picture from $\pi$ to $\pi'$,
and denote by $\omega,\omega'$ the respective labellings. Then, the
first condition on a picture means that $\omega'\circ \phi$ is an
increasing bijection of $(P,<_1)$ into $\{1,...,n\}$; hence
$\omega\circ \phi^{-1}\circ\omega'^{-1}$ is a linear extension of
$(P,<_1)$. Similarly, the second condition means that $\omega'\circ
\phi\circ\omega^{-1}$ is a linear extension of $(P',<_1)$.
Thus the lemma follows.
 \hfill $\Box$

\noindent{\bf Proof of theorem.} The fact that $L$ preserves the
scalar product is immediate from the lemma.

It remains to show that $L$ is a homomorphism for the internal
product. Let $\pi, \pi'$ be two special double posets with
underlying sets $E,F$. We may assume that $E=F=\{1,...,n\}$, such
that their second order $<_2$ is the natural order on $\{1,...,n\}$.
Then, since the labellings of $E$ and $F$ are the identity mappings, a linear extension of $\pi$ (resp. $\pi'$) is a
permutation $\alpha$ (resp. $\beta$) in $S_n$ such that
$\alpha^{-1}$ (resp. $\beta^{-1}$) is increasing from $(E,<_1)$
(resp. $(F,<_1)$) into $\{1,...,n\}$.

Now, let $\phi$ be increasing $(E,<_1)\rightarrow
(F,<_2)=\{1,...,n\}$. We construct the double poset
$\Pi=E\times_\phi F$ as in Section 2.2. Since we identify isomorphic
double posets, we may take $F$ as underlying set, with the first
order $<_1$ of $F$ as first order of $\Pi$, and with second order
defined by the labelling $\phi^{-1}: F\rightarrow E=\{1,...,n\}$.
Then a linear extension $\sigma$ of $\Pi$ is a permutation $\sigma$
such that $\sigma{-1}\circ\phi^{-1}$ is increasing
$(F,<_1)\rightarrow \{1,...,n\}$. Define $\alpha=\phi^{-1}$ and
$\beta=\phi\circ\sigma$. Then $\alpha^{-1}=\phi$ (resp.
$\beta^{-1}=\sigma^{-1}\circ\phi^{-1}$) is increasing
$(E,<_1)\rightarrow \{1,...,n\}$ (resp. $(F,<_1)\rightarrow
\{1,...,n\}$), and therefore $\alpha$ and $\beta$ are linear
extensions of $\pi$ and $\pi'$ with $\alpha\circ\beta=\sigma$.

Conversely let $\alpha$ and $\beta$ be linear extensions of $\pi$
and $\pi'$. Put $\sigma=\alpha\circ\beta$ and $\phi=\alpha^{-1}$.
Then $\phi$ is increasing $(E,<_1)\rightarrow (F,<_2)=\{1,...,n\}$
and $\sigma^{-1}\circ\phi^{-1}=\beta^{-1}$ is increasing
$(F,<_1)\rightarrow \{1,...,n\}$. Hence $\sigma$ is a linear
extension of $\Pi$.

All this implies that $L(\pi)L(\pi')$, which is the sum of all
$\alpha\circ\beta$'s, is equal to $L(\pi\circ\pi')$, which is equal
to the sum of all $\sigma$'s, for all possible $\phi$'s .
 \hfill $\Box$
\medskip

It is easy to prove also that the submodule of $\BZ\mathbf {DS}$
spanned by the {\it naturally labelled} special double posets, that
is, those whose second order is a linear extension of the first
order, is a sub-bialgebra of $\BZ\mathbf {DS}$ (the class of naturally
labelled posets is closed under taking disjoint unions and ideals).
It may be possible that one could compute the antipode of this subalgebra,
and that of $\BZ\mathbf {DS}$, by extending the techniques of Aguiar and Sottile
\cite{ASo}, who computed the antipode of $\BZ S$ using the weak
order on the symmetric group, and a recursive method in posets, as
\cite{Gess} proof of Th.1 and \cite{BlSc} Lemma 4.11.

Recall that for a permutation $\sigma\in S_n$, its {\it descent
composition} $C(\sigma)$ is the composition of $n$ equal to
$(c_1,...,c_k)$, if $\sigma$ viewed as a word has $k$ consecutive
ascending runs of length $c_1...,c_k$ and $k$ is minimum. For example
$C(51247836)=(1,5,2)$, the ascending runs being $5,12478, 36$.
Recall from \cite{Gess} p.291 the definition of the fundamental
quasi-symmetric function $F_C$, for any composition $C$ (see also
\cite{Stan1} 7.19 where it is denoted $L_\alpha$). Then it follows
from \cite{MaRe} Th.3.3 that the linear function $F: \BZ
S\rightarrow \mathbf {QSym}$ defined by $\sigma\mapsto
F_{C(\sigma)}$ is a homomorphism of bialgebras. Recall that the
bialgebra homomorphism $\Gamma: \BZ{\mathbf D}\rightarrow \mathbf
{QSym}$ has been defined in Section 3. Then the following result is
merely a reformulation of a result of Stanley (see \cite {Gess} Th.1
and Eq.(1) page 291, or \cite{Stan1} Cor.7.19.5).

\begin{corollary}
The mapping $F\circ L$ is equal to $\Gamma$ restricted to $\BZ\mathbf {DS}$.
\end{corollary}

\subsection{Littlewood-Richardson rule}

A {\it lattice permutation} (or Yamanouchi word) is a word on the
symbols in ${\mathbb P} = \{1,2, 3, ...\}$ such that, for any  $i$,
in each left factor, the number of $i$'s is not less than the number
of $i+1$'s. For example, $11122132$ is such a lattice permutation.
Given a word $w$ on the symbols $1,2, 3, ...$, with $k$ the greatest
symbol appearing in it, we call {\it complement} of $w$, the word
obtained from $w$ by exchanging $1$ and $k$ in $w$ , then $2$ and
$k-1$, and so on. For the word of the above example,  its complement
is therefore $33322312$. The {\it weight} of a word $w$ is  the
partition $\nu=1^{n_1}2^{n_2}...$, where $n_i$ is the number of
$i$'s  in $w$. For the word above, it is the partition $1^42^33^1$.

Given a special double poset $\pi=(E,<_1,<_2)$ of cardinality $n$
with labelling $\omega$, and a word $a_1a_2...a_n$ of length $n$
over a totally ordered alphabet $A$, we say that $w$  {\it  fits
into} $\pi$ (we take the terminology from \cite{GR}) if the function $E\rightarrow  A$ defined by $e\mapsto
a_{\omega(e)}$  is a $\pi$-partition. In other words, given a
function $f: E\rightarrow A$, call {\it reading word} of  $f$ the
word $f(\omega(1))...f(\omega(n))$. Then $f$ is a $\pi$-partition if
and only if its reading word fits into $\pi$. Note the case where
the word is a permutation $\tau$: we have that $\tau$ fits into
$\pi$ if and only if the word $\tau(1)...\tau(n)$ fits into $\pi$. This means
that the mapping $e\mapsto\tau(\omega(e))$ is a $\pi$-partition,
that is, since it is a bijection, is increasing  $(E,<_1)\rightarrow
\{1,...,n\}$.

Given a partition $\nu$  of $n$, we define (as in \cite{Gess}) a special double poset
$\pi_\nu=(E_\nu, <_1,<_2)$  where $E_\nu$ is the Ferrers diagram of
$\nu$, where $<_1$ is the  order induced on $E_\nu$ by the natural
partial order of  $\BN\times\BN$, and where $<_2$ is given on the
elements of  $E_\nu$ by $(x,y)<_2(x',y')$ if and only if either $y>y'$, or
$y=y'$  and $x<x'$. Recall that there is a
well-known bijection between standard Young tableaux of shape $\nu$
and lattice permutations of weight $\nu$, see \cite{Stan1}
Prop.7.10.3 (d). In this bijection, the shape of the tableau is
equal to the weight of the lattice permutation.

If $\pi$ is a special double poset, we denote by $\tilde \pi$ the
special double poset obtained by replacing the two orders of $\pi$
by their opposite. Clearly, a permutation $\sigma$ fits into $\pi$
(assumed to be special) if and only if $w_0\circ \sigma \circ w_0$
fits into $\tilde \pi$ ($w_0$ is defined in Lemma 3.3 below).

\begin{theorem} Let $\pi$ be a special double poset and $\nu$ be some partition.
Then, the scalar product  $(\pi,\pi_\nu)$, that is, the number of
pictures from $\pi$ to $\pi_\nu$, is equal to:

(i) the number of lattice permutations  of weight $\nu$ whose complements fit into $\pi$;

(ii) the number of lattice permutations of weight $\nu$ whose mirror images fit into $\tilde \pi$.
\end{theorem}

Note that part (ii) of this is the classical formulation of the
Littlewood-Richardson rule (see \cite{Ma} (9.2) or \cite{Stan1}
Th.A1.3.3), once one realizes that a skew Schur function indexed by
a skew shape is equal to the skew Schur function obtained by
rotating by 180 degrees that shape (cf. \cite{BlSc} Chapter 11
p.109-110).

We first need the following:
\begin{lemma}
Let $\pi$ be a special double poset.  A permutation $\sigma$ is a
linear extension of $\pi$ (with respect to $<_1$) if and only if its
inverse fits into $\pi$.
\end{lemma}

\noindent{\bf Proof.} We may assume that $\pi=(E,<_1,<_2)$ with
$E=\{1,...,n\}$ and $<_2$ the natural order of $E$. Then $\omega $
is the identity and therefore a permutation $\sigma$ is a linear
extension of  $\pi$ if and only if $\sigma^{-1}$ is increasing
$(E,<_1)\rightarrow \{1,...,n\}$. On the other hand, $\tau$ fits
into $\pi$ if and only if it is increasing $(E,<_1)\rightarrow
\{1,...,n\}$, as noted previously. \hfill $\Box$

\begin{prop}
Let $\pi$ be a special double poset.  A word $w=a_1...a_n$ fits into
$\pi$ if and only if its standard permutation does.
\end{prop}

This result is equivalent to a result of Stanley, see \cite{Gess}
Th.1 or \cite{Stan1} Th.7.19.14.

A standard Young tableau of shape $\nu$ is the same thing as a
$\pi_\nu$-partition which is a bijection from the Ferrers diagram of
$\nu$ onto $\{1,...,n\}$. Thus we can speak of the reading word of a
tableau, which is a classical notion. We denote it $read(T)$. We
consider also the {\em mirror reading word} of $T$, which is the mirror image
of $read(T)$. The following result must be well-known, but we give a
proof for the convenience of the reader. Part (i) is proved in
\cite{BlSc} page 109.

\begin{prop} Let $T$ be a standard Young tableau and $w$ the associated lattice permutation.

(i) Let $u$ be the complement of $w$. Then the reading word of $T$
is equal to the inverse of the standard permutation of $u$.

(ii) Let $v$ be the mirror image of $w$.  Then the mirror reading
word of $T$ is equal to the complement of the inverse of the
standard permutation of $v$.
\end{prop}

Consider for instance the Young tableau
$$
\begin{array}{cccc}
5 \\
3&9 \\
2&6&10&11 \\
1&4&7&8
\end{array}
$$
Then its lattice permutation is $w=1\ 2\ 3\ 1\ 4\ 2\ 1\ 1\ 3\ 2\ 2$.
The complement of the latter is $u=4\ 3\ 2\ 4\ 1\ 3\ 4\ 4\ 2\ 3\ 3$.
Standardizing this latter word, we obtain the permutation $8\ 4\ 2\
9\ 1\ 5\ 10\ 11\ 3\ 6\ 7$, whose inverse is $5\ 3\ 9\ 2\ 6\ 10\ 11\
1\ 4\ 6\ 7$, which is indeed the reading word of the given tableau,
obtained by concatenating its rows, beginning with the last row.
This illustrates (i). For (ii), the mirror image of $w$ is $v=\ 2\
2\ 3\ 1\ 1\ 2\ 4\ 1\ 3\ 2\ 1$. The standard permutation of $v$ is $\
5 \ 6 \ 9 \ 1 \ 2 \ 7 \ 11 \ 3 \ 10 \ 8 \ 4$. The inverse of this
permutation is $\ 4\ 5\ 8\ 11\ 1\ 2\ 6 \ 10 \ 3 \ 9\ 7$. Finally,
the complement of this permutation is $\ 8\ 7\ 4\ 1\ 11\ 10\ 6\ 2\
9\ 3\ 5$, which is indeeed the mirror reading word of $T$.

We need a lemma. For this, we use a variant of the reading word of a
tableau. Call {\it row word} of a standard Young tableau $T$ the
permutation, in word form,  obtained by reading in increasing order
the first row of $T$, then the second, and so on. Denote it by
$row(T)$. For example, the row word of the previous example is
$$row(T)= 1 \ 4\ 7\ 8\ 2\ 6\ 10\ 11\ 3\ 9\ 5.$$

\begin{lemma} Let $w$ be a word on the alphabet $\mathbb P$
of weight equal to the partition $\nu=(\nu_1 > ... > \nu_k > 0)$ of
$n$. Let $w_0= n\ n-1\ \ldots 2\ 1$ be the longest element in the
group $S_n$ and $\gamma$ the longest element in the Young subgroup
$S_{\nu_1}\times...\times S_{\nu_k}$.

(i) Let $u$ be the complement of $w$. Then $st(u)=
w_0\circ\gamma \circ st(w)$.

(ii) Let $v$ be the mirror image of $w$. Then $st(v)=
\gamma \circ st(w) \circ w_0$.

(iii) Suppose that $w$ is a  lattice permutation and let $T$ be the
tableau of shape $\nu$ corresponding to $w$. Then $st(w)$ is the
inverse of $row(T)$.
\end{lemma}

\noindent{\bf Proof.} (i) Note that $\gamma$ as word is equal to
$$\nu_1 \ldots 1 \ \ (\nu_1+\nu_2) \ldots (\nu_1+1)\ \ \ldots \ \ n \ldots
(\nu_1+\ldots+\nu_{k-1}+1).$$ Hence
$$w_0\circ\gamma = (\nu_2+\ldots+\nu_k+1)\ldots n \ \
(\nu_3+\ldots+\nu_k+1)\ldots (\nu_2+\ldots+\nu_k) \ \ \ldots \ \ 1
\ldots \nu_k.$$ The proof of (i) then follows by inspection.

(ii) Likewise, one proves that $st(v)\circ w_0=\gamma \circ st(w)$ by inspection.

(iii) Let $I_1, \ldots, I_k$ denote the successive intervals of
$\{1,\ldots,n\}$ of cardinality $\nu_1, \ldots, \nu_k$. Let
$L_1,\ldots, L_k$ be the set of elements in the successive rows of
$T$. If $I$, $L$ are two subsets of equal cardinality of $\mathbb
P$, we denote by $I\nearrow L$ the unique increasing bijection from
$I$ into $L$. We denote also $f_1\cup \ldots \cup f_k$ the function
which restricts to $f_i$ on its domains, assuming the domains are
disjoint. Then $row(T)$ is the permutation $\cup_{j=1,\ldots,k}
I_j\nearrow L_j$, and its inverse is therefore $\cup_{j=1,\ldots,k}
L_j\nearrow I_j$. Moreover, the word $w$ is defined by the following
condition: for each position
$p\in\{1,\ldots,n\}=\cup_{j=1,\ldots,k}L_j$, the $p$-th letter of
$w$ is $j$ if and only if  $p\in L_j$. Recall that $st(w)$, viewed
as word, is obtained by giving the numbers $1,...,n$ to the
positions of the letters in $w$, starting with the $1$'s from left
to right, then $2$'s, and so on. Therefore
$st(w)=\cup_{j=1,\ldots,k} L_j\nearrow I_j.$ \hfill$\Box$
\medskip

\noindent{\bf Proof of proposition.} (i) We have to prove that
$read(T)=st(u)^{-1}$. We know by Lemma 5.2.(i) and (iii) that
$st(u)=w_0 \circ\gamma \circ st(w)$ and $st(w)^{-1}=row(T)$. Clearly
$read(T)=row(T)\circ \delta$, where $\delta$ is the permutation
$$(\nu_1+\ldots+\nu_{k-1}+1)\ldots  n \ \ (\nu_1+\ldots+\nu_{k-2}
+1)\ldots(\nu_1+\ldots+\nu_{k-1}) \ \ \ldots \ \ 1\ldots \nu_1.$$
Now $\delta=\gamma \circ w_0$. Therefore, $read(T)=row(T) \circ
\gamma \circ w_0=st(w)^{-1}\circ\gamma\circ w_0=st(u)^{-1}$, since
$w_0$ and $\gamma$ are involutions.

(ii) We have to show that $read(T)\circ w_0=w_0 \circ st(v)^{-1}$.
We know by Lemma 5.2.(ii) that $st(v)= \gamma \circ st(w) \circ
w_0$, hence $st(w)=\gamma \circ st(v) \circ w_0$. Using what we have
done in(i), we have therefore $$read(T)= st(w)^{-1} \circ \gamma
\circ w_0 = w_0 \circ st(v)^{-1} \circ\gamma \circ\gamma \circ
w_0=w_0 \circ st(v)^{-1}\circ w_0. $$ \hfill $\Box$
\medskip

\noindent{\bf Proof of theorem.} (i) By Lemma 3.1 and Lemma 3.2, the
indicated scalar product is equal to the number of permutations
$\sigma$ which fit into $\pi$ and whose inverse fits into $\pi_\nu$.
Let $H$ denote the set of complements of lattice permutation of
weight $\nu$. By Proposition 3.2.(i), the mapping $H\rightarrow
RWSYT_\nu$, $w\mapsto st(w)^{-1}$ is a bijection, where we denote by
$RWSYT_\nu$ the set of reading words of standard Young tableaux of
shape $\nu$. By Proposition 3.1, part (i) of the theorem follows.

(ii) Let $K$ denote the set of mirror images of lattice permutations of weight $\nu$. By Prop. 3.2.(ii), the mapping $K\rightarrow RWSYT_\nu$, $v\mapsto \sigma=w_0\circ st(v)^{-1}\circ w_0$ is a bijection. Moreover, $\sigma^{-1}=w_0\circ st(v)\circ w_0$ fits into $\tilde \pi$ if and only if $st(v)$ fits into $\pi$, that is, by Prop. 3.1, if and only if $v$ fits into $\pi$.
\hfill $\Box$

\end{document}